\documentclass[a4paper,notitlepage, 8pt,reqno]{article}

  \usepackage{mathtext}
 \usepackage[cp1251]{inputenc}
  \usepackage[T2A]{fontenc}
 \usepackage[english]{babel}
 \usepackage[all,arc,poly,2cell,curve,arrow,tips]{xy}

  \usepackage{enumerate, eucal, amsthm,amsmath, amssymb}

  \allowdisplaybreaks[4]

 \theoremstyle{definition}
 \newtheorem{defn}{Defenition}%[section]

 \theoremstyle{plain}
 \newtheorem{thm}{Theorem}
 \newtheorem*{thm*}{Theorem}
 \newtheorem{prop}{Proposition}
 \newtheorem{lem}{Lemma}

 \theoremstyle{remark}
  \newtheorem{remark}[defn]{Remark}

 \renewcommand{\abstractname}{}

  \newcounter{ab}
\title{
An estimate of the number of apparent singularities in the
Riemann-Hilbert problem on a compact Riemann surface.
 }

 \author{D. V.  Artamonov}

\begin{document}
 \maketitle

\renewcommand{\abstractname}{}

\begin{abstract}
In the paper we give an upper estimate of the number of apparent
singularities that are sufficient for construction of a system of
regular linear differential equations on a Riemann surface with
given fuchsian singularities and monodromy.
\end{abstract}

 Let $M$ be a compact Riemann surface of genus $g>0$, let $a_1,...,a_n$, $a_i\neq a_j$ be points on this surface.  The fundamental group $\pi_1(M\setminus\{a_1,...,a_n\};x_0)$ of the surface $M$ with punctures $a_1,...,a_n$
  is described as follows. It has generators $\gamma_1,...,\gamma_n,f_1,h_1,...,f_{g},h_{g}$  and one relation:
   $\gamma_1...\gamma_nh_1f_1h_{1}^{-1}f_1^{-1}...h_{g}f_{g}h_{g}^{-1}f_{g}^{-1}=1
  $.  The generators   $f_1,h_1,...,f_{g},h_{g}$ are represented
 by the same loops, that give generators of  $\pi_1(M;x_0)$, and generators $\gamma_i$
  are represented by loops starting from $x_0$ than going around the point $a_i$ and than returning to $x_0$
  (but in such a way that these loops are contractible in
  $M$).

 Let us be given a representation $\chi:\pi_1(M\setminus\{a_1,...,a_n\};x_0)\rightarrow GL_p(\mathbb{C})$,  that is a collection of matrices
 $G_1,...,G_n,F_1,H_1,...,F_{g},H_{g}$, satisfying the relation
  $G_1...G_nH_1F_1H_{1}^{-1}F_1^{-1}...H_{g}F_{g}H_{g}^{-1}F_{g}^{-1}=E
  $.
 Does there exist a system of  $p$ linear differential equations on $M$ with
  fucshian singularities in $a_1,...,a_n$, such that $\xi$ - is its monodromy representation?

This system is determined by a differential form with poles of the
first order in $a_1,...,a_n$. If $g>0$,  the dimension of the space
of such forms is lower than the dimension of the space of
representations of fundamental group.  So, when $g>0$,  the answer
is typically negative.

In turns out, that every representation can be realized as a
monodromy representation  in a class of regular systems with
additional singularities $b_1,...,b_N$ with a trivial monodromy
corresponding to a bypass around each of them. That means the
following. For every representation
$\chi:\pi_1(M\setminus\{a_1,...,a_n\};x_0)\rightarrow
GL_p(\mathbb{C})$ there exists a system of $p$ differential
equations with regular  singularities $a_1,...,a_n,b_1,...,b_N$ with
a trivial monodromy corresponding to a bypass around apparent
singularities $b_1,...,b_N$. A singularity is regular if every
solution of a system as a power-like behavior, when we approach to
the singularity inside some cone (that is  we are not allowed to
rotate around the singularity).

Estimates of the number $N$ of additional singularities were given
by different authors. In \cite{O} such an estimate was given for the
problem of construction of a fucshian differential equation.  But in
that article some strong restrictions on the monodromy
representation were  suggested and the resulting estimate depends
not only on $p$ and $g$, but also on the number of singularities
$n$.  In articles \cite{Y},\cite{B} the case $p=2$ was considered.
But  some additional restrictions on the monodromy representation
were made.

Apparent singularities are always  regular. In the article \cite{B}
in the situation under consideration the orders of the poles of the
coefficients of the system in these apparent singularities are
estimated (in the given singularities $a_i$ the system has fuchsian
singularities). We shall not do this in the present article.

Also it is known, that if $g=0$, than it is sufficient to add one
(or none) additional singularity. That's why we shall suggest that
$g>0$.

We shall prove the following main theorem:

\begin{thm} Every representation of $\pi_1(M\setminus\{a_1,...,a_n\};x_0)$ in $GL_p(\mathbb{C})$ can be realized as a monodromy representation of a regular system with less or equal than
 $2pg-g+1$ additional singularities with trivial monodromy corresponding to a bypass around every additional singularity.
\end{thm}

\proof
 We are given a representation
$\chi:\pi_1(M\setminus\{a_1,...,a_n\};x_0)\rightarrow
GL_p(\mathbb{C})$. For this representation we construct a
holomorphic vector bundle $E$ with a connection $\nabla$, such that
$a_1,...,a_n$ are fucshian singularizes of the connection and the
monodromy of that connection is $\chi$ (\cite{B1}).

{\bf The first step} of the proof is to show the connection between
meromorphic trivialization of the bundle $E$ and constructions of a
system a differential equations with monodromy $\xi$ with some
additional singularities.

 It is known that on a Riemann sphere every holomorphic vector bundle splits into a direct sum of linear bundles (the Birkhoff-Grothendieck theorem).
 For Riemann surfaces of higher genus that is not true. Nevertheless, the meromorphic triviality of every bundle takes place. That means the following:

  \begin{thm*}(\cite{F})  For every $p$-dimensional  holomorphic vector bundle there exist meromorphic sections  $\psi_1,...,\psi_p$ such that they form a bases in
  stalks over  all points of the base except a finite number of points
  of the
  base.
  \end{thm*}

 Among the  points mentioned in the theorem there are poles of  $\psi_1,...,\psi_p$  and points $z_1,...,z_t$, where
   all the sections $\psi_1,...,\psi_p$ are finite but they are
  linearly dependent.

%\begin{prop}
A meromorphic section $s$ of a bundle with meromorphic
trivialization  $\psi_1,...,\psi_p$ is uniquely  defined by a
collection of meromorphic functions $(\alpha_1,...,\alpha_p)$ by the
formula $s=\alpha_1\psi_1+...+\alpha_p\psi_p$.
%\end{prop}
%\proof

%Of course the formula $s=\alpha_1\psi_1+...+\alpha_p\psi_p$ defines
%a meromorphic section.  We have to understand, why every section is
%given by such a formula.

%We are given a section  $s$.  Outside of poles and points
%$z_1,...,z_t$ the
% sections $\psi_1,...,\psi_p$  form a bases in stalks, so  $s$ on
% this set
% equals to $\alpha_1\psi_1+...+\alpha_p\psi_p$, where $\alpha_i$ are meromorphic functions.  We have only to
% understand, what kind of singularities these functions  have in poles of sections $\psi_1,...,\psi_p$  and in the points $z_1,...,z_t$.
% In some neighborhood  of these points the section is given by the formula  $s=\beta^i_1v_1^i+...\beta_p^iv_p^i$, where
%$v_1^i,...,v_p^i$   are   holomorphic base sections defined locally
%in a neighborhood of  a considered point. As $s$  is meromorphic,
%functions $\beta^i_1,...,\beta^i_p$ are meromorphic.
% Let's denote as  $\beta^i$ the the vector  $(\beta^i_1,...,\beta^i_p)$, and as  $\alpha$ the vector  $(\alpha_1,...,\alpha_p)$.
%Then $s=\alpha\psi=\alpha W^iv^i=\beta^i v^i$, where $W^i$ is a
%transition matrix form  the bases $(v_1^i,...,v_p^i)$ to the bases
%$(\psi_1,...,\psi_p)$. We get that $\alpha=\beta^i (W^{i})^{-1}$.
%But the matrix $(W^{i})^{-1}$ is meromorphic so the vector $\alpha$
%is also meromorphic.

%\endproof

Now we are going to explain the relation between the meromorphic
trivializations for bundles and an estimate of apparent
singularities. Given a representation $\chi$, we have constructed a
bundle $E$ and a connection $\nabla$ such that the  monodromy of
$\nabla$ is $\chi$. Let's find some meromorphic trivialization for
$E$.

\begin{prop}
Let $(E,\nabla)$ be a vector bundle with a fuchsian connection
constructed from singularities $a_1,...,a_n$ and  a monodromy
$\chi$. Let $\psi_1,...,\psi_p$ be meromorphic a trivialization of
$E$. Then there exist a  regular system of linear equations with the
monodromy $\chi$ and a set of additional singularities contained in
the following set $S$: poles of trivializing sections $\psi_i$,
points $z_1,...,z_t$, where they are linearly dependent. If
$a_i\notin S$, then $a_i$ is a fuchsian singularity for the system.
\end{prop}

\proof  The sections $\psi_1,...,\psi_p$  form a bases in stalks of
$E$ over points from  $M\setminus S$.

 A meromorphic section is defined by a vector-function $(\alpha_1,...,\alpha_p)$.
The condition that this section $(\alpha_1,...,\alpha_p)$ is
horizontal is written as a system of linear differential  equations
which outside $S\cup\{a_1,...,a_n\}$ is nonsingular. The monodromy
of that system coincides with the monodromy of the connection. So we
get a system with needed monodromy and it's additional singularities
contain in $S$.

\endproof

{\bf The second step}. We formulate our problem in the following
way. Given a $p$-vector  bundle $E$, we must find a meromorphic
trivialization with as least as possible  number of points that are
poles of trivializing sections or points, where these sections are
finite, but linearly dependent.

\begin{remark} Such an estimate is known, if $E$ is a  stable bundle. From the technique of Turin
parameters (their definition can be found in \cite{KN}, see also the
original work (\cite{Tu}) it follows, that for a stable bundle $pg$
 additional singularities are enough. But not every
representation can be realized as a representation of a connection
in a stable (or even in a semistable) bundle (see \cite{B1})
\end{remark}

Before formulate our plan, we shall prove the following:
 \begin{prop}
On a Riemann surface $M$ of genus  $g$ for every point $P\in M$
every linear bundle
 $L$ has a meromorphic section $s$ with the following property. On
$M\setminus P$ it has less or equal than $g$  zeroes and no poles.
Also  $(s)+(-degL+g)P\geq 0$, where $(s)$ is a divisor of zeroes and
poles of section $s$.
 \end{prop}
The last statement can be reformulated as follows:  in $P$ the
section $s$ has a pole of order less or equal than
 $(-degL+g)$.
 \proof

 If we multiply our bundle by the bundle $\xi(P)^{k}$, we get a bundle
   $L'=L\otimes \xi(P)^{k}$  with degree
  $degL'=degL+k$.  Put $k=-degL+g$.
By Riemann-Roch theorem  we get an inequality
$dimH^0(X,\mathcal{O}(L'))\geq (1-g)+deg L'=1>0$.  So the bundle
$L'$  has a holomorphic section $s'$.

As $s'$ is holomorphic, it has no poles, but it has $degL'=g$ zeroes
with multiplicities and less or equal than $g$ points were it takes
zero values.  From this section of the bundle $L'=L\otimes
\xi(P)^{k}$ we construct a section  $s$ of  the bundle $L$ in the
following way. Every section of  $L'=L\otimes \xi(P)^{k}$ has the
form $s'=ss_0$, where $s$ is a section of $L$, and $s_0$ is a
section of $\xi(P)^{k}$. From this we get that  $\frac{s'}{s_0}$  is
a section of $L$.  Let's take as $s_0$ a canonical section of
$\xi(P)^{k}$, which has a zero of order $k$  (or a pole if $k<0$) in
$P$.

The section $s'$  has less or equal than $g$  zeroes.  From explicit
formulas we see that on $M\setminus P$  the section $s$  also has
less or equal  than $g$ zeroes.

 The divisor $(s)$  of zeroes and poles of $s$ equals to $(s')-(-degL+g)P$. From here and the fact that
$(s')\geq 0$, the second statement follows.

\endproof

Now let's formulate, what we shall do.
 Let's fix some point $P\in M$. Using the
previous proposition, we shall suggest that the number of zeroes of
every trivializing section $\psi_i,i=1,...,p$  is less or equal than
$g$ and all these section have probably only a single pole in $P$.

 So the number of poles is minimized and we have to minimize only
 the number of points. We shall inductively construct new
 trivializing sections $\psi''_1,...,\psi''_p$. The will have the
 following properties.

 \begin{enumerate}
 \item $<\psi''_1,...,\psi''_k>=<\psi_1,...,\psi_k>$
 \item All sections  $\psi''_k$ can have pole only in the point
 $P$.
 \item The number of points, where all these sections
 $\psi''_1,...,\psi''_k$ are finite, but linearly dependent is less
 or equal than the some number $N(p,g)$, which depends only on $p$
 and $g$.

 \end{enumerate}

Also we shall calculate the number $N(p,g)$,  it turns out, that it
equals to $2pg-g$.

{\bf In the third step} we shall explain, how we construct
$\psi''_1,...,\psi''_p$. We shall construct them by induction by
$p$. Let $\psi''_1=\psi_1$, so $N(1,g)=g$. Suggest that an estimate
for $N(p-1,g)$ is done, let's do it for $N(p,g)$.

Let's consider subbundle $E'\subset E$ formed by sections
$\psi_1,...,\psi_{p-1}$. We have already constructed meromorphic
sections $\psi''_1,...,\psi''_{p-1}$
 of  $E'$,  which give meromorphic trivialization, so that the number of points $z_1,...,z_m$,
  where these sections are finite but linearly dependent is less or equal than $N(p-1,g)$.  All these sections, may be, have a single pole in $P$.

If we add to $\psi''_1,...,\psi''_{p-1}$ the section  $\psi_p$ - we
get a meromorphic trivialization of $E$. Where all these sections
are linearly dependent?  First of all, in points  $z_1,...,z_m$,
where sections $\psi''_1,...,\psi''_{p-1}$ are dependent. Then in
points $z_{m+1},...,z_t$, where  $\psi''_1,...,\psi''_{p-1}$ are
independent, but  $\psi''_1,...,\psi''_{p-1}, \psi_p$  are
dependant.

 Let's take a point $z_i,\,\,\,i>m$. As sections
 $\psi''_1,...,\psi''_{p-1}$  are independent, than as a bases holomorphic sections in a neighborhood of  this point we can take sections  $
 \psi''_1,...,\psi''_{p-1},v^i$, where $v^i$ is some holomorphic section, defined in a neighborhood of $z_i$.

In that neighborhood we have equality
$\psi_p=\alpha_1^i\psi''_1+...+\alpha^i_{p-1}\psi''_{p-1}+\alpha^i_{p}v^i$
for some holomorphic functions $\alpha_1^i,...,\alpha^i_{p}$,
defined in a neighborhood of $z_i$.

 As sections $\psi''_1,...,\psi''_{p-1},\psi_p$  are linearly independent in  $z_i$,
than in that point $\alpha^i_p$  has a zero with multiplicity
$d_i>0$,  that is $\alpha^i_p=z^{d_i}{\alpha'}^i_p$, where  $z$  is
a local coordinate in a neighborhood of $z_i$, such that $z(z_i)=0$.
Here ${\alpha'}^i_p(z_i)\neq 0$.

Let prove:
\begin{lem}
There exist meromorphic functions
$\widetilde{\alpha_1},...,\widetilde{\alpha_{p-1}}$ on $M$, such
that they have a single pole in $P$ and for every $i={m+1},...,t$,
except less or equal than $g$ indices,  in $z_i$  the function
$\widetilde{\alpha_j}-\alpha_j^i$, $j=1,...,p-1$ has  a zero of
order at least
 $d_i+1$.
\end{lem}
In other words in these points the first $d_i+1$ members of a Teilor
series of a global function $\widetilde{\alpha_j}$ and of a local
function $\alpha_j^i$  must coincide.

\begin{defn} Those points $z_{m+1},...,z_t$,  in which the function
$\widetilde{\alpha_j}-\alpha_j^i$ has zero of needed order will be
called unexceptional.
\end{defn}
So there are  less or equal  than $g$  exceptional points.

 \proof

Let $f$ be a function which has probably a pole in  $P$  and zeroes
of the first order in points $z_{m+1},...,z_t$, except less or equal
than $g$ points.  It is constructed as follows:  we take a
meromorphic section $s$ of a bundle
$\xi(z_{m+1})^{-1}\otimes...\otimes\xi^{-1}(z_t)$, which has, maybe,
a single pole in $P$  and less or equal  than $g$ zeroes. This
bundle also has a canonical section $s^0$ with first order poles in
points $z_{m+1},...,z_t$ and no zeroes. Than $f=\frac{s}{s^0}$ is a
meromorphic function.

 If $z_i$ is not a zero of  $s$, than in  $z_i$ the function $f$ has a zero of the first order.
  If in $z_i$  the section $s$  has zero, than $f$ has in this point a zero of a higher order. This point is
  exceptional, the number of such points is less or equal than the number of zeroes of $s$,  which is less or equal than than $g$. The pole of $f$ is placed at the  point   $P$.

 Let's now construct  functions $f_{m+1},...,f_t$, which have poles in $P$, and such that $f_i$  takes nonzero value in $z_i$  and zero  values in other points
$z_{m+1},...,\widehat{z_i},...,z_t$.

To construct such functions we take a bundle
$\xi(z_{t+1})^{-1}\otimes...\otimes\widehat{\xi^{-1}(z_i)}\otimes...\otimes\xi^{-1}(z_t)$.
It has a section $s$ with a pole in $P$.  Also it has a canonical
section $s^0$ with poles of the first order in
$z_{m+1},...,\widehat{z_i},...,z_t$  and no zeroes.  Let's define a
meromorphic function $\widetilde{f}_i=\frac{s}{s^0}$
  That function has a single pole in $P$  and zeroes in
$z_{t+1},...,\widehat{z_i},...,z_m$. If it turns out that
$\widetilde{f}_i(z_i)\neq 0$,  than we put $\widetilde{f}_i=f_i$.

If $\widetilde{f}_i(z_i)= 0$, than we do the following. Let the
order of zero of a function $\widetilde{f}_i$ in $z_i$ be equal to
$k$. There exists a meromorphic function $h$,  which has a single
pole in $z_i$ and it's order equals to $r$.  Put
$f_i=(\widetilde{f}_i)^rh^k$.  This function takes nonzero value in
$z_i$, and has a single pole in  $P$,  where $\widetilde{f}_i$,
takes zero values in points $z_{m+1},...,\widehat{z_i},...,z_t$,
where $\widetilde{f}_i$ has zeroes.

Let $d$ be $max_i\{d_i+1\}$, denote as $g_i$  the function $f_i^d$.
This function takes nonzero value in $z_i$  and has zeroes in
$z_{m+1},...,\widehat{z_i},...,z_t$  of orders at least $d$. The
poles of these functions are in the point  $P$.

Now we are ready to construct $\widetilde{\alpha_j}$. For
unexceptional point $z_i$ we have a function $g_i$,  which takes
nonzero value in $z_i$, and a function $fg_i$, which has a zero of
the first order in  $z_i$. In other points
$z_{m+1},...,\widehat{z_{i}},...,z_{t}$ (in exceptional also)  it
has zeroes of orders not less then $d$. There exist a polynomial
without constant term $p^i_j$, such that $p^i_j(fg_i,g_i)$ and
$\alpha_j^i$ have the same first $d_i+1$ members of Teilor series in
point $z_i$. As in other points $z_{m+1},...,\widehat{z_i},...,z_t$
all functions $g_i,fg_i$ have zeroes of order at least
$d=max\{d_i+1\}$, than the function $p^i_j(fg_i,g_i)$ has zero first
$d$ members of Teilor series in these points.

Thus $\widetilde{\alpha_j}=\sum_{i}p^{i}_j(fg_{i},g_i)$, where the
summation is taken over all unexceptional points, is a needed
function.
\endproof

Let's return to the number $N(p,g)$. First of all, instead of
$\psi_p$ we take a section
$\psi'_p=\psi_p-\widetilde{\alpha}_1\psi''_1-...-\widetilde{\alpha}_{p-1}\psi''_{p-1}$.
All the functions $\widetilde{\alpha}_i$ and sections $\psi''_1,
...,\psi''_{p-1},\psi_p$ can have a pole only in $P$. That is why
$\psi'_p$ has, maybe,  a pole only in this point $P$.

As the transformation $\psi_p\mapsto \psi'_p$ is invertible
everywhere except $P$, that in the points where sections
$\psi''_1,...,\psi''_{p-1},\psi_p$ are independent, the sections
$\psi''_1,...,\psi''_{p-1},\psi'_p$ are also independent. So
$\psi''_1,...,\psi''_{p-1},\psi'_p$ are dependent only in points
$z_1,...,z_t$.

By induction, we have minimized the number of points, where
$\psi''_1,...,\psi''_{p-1}$  are dependent and we have estimated it
by  a number $N(p-1,g)$. In the neighborhood of a point
$z_i\in\{z_{m+1},...,z_t\}$ there were equalities
$\psi_p=\alpha_1^i\psi''_1+...+\alpha^i_{p-1}\psi''_{p-1}+\alpha^i_pv^i$.
They give equalities
$\psi'_p=(\alpha^i_1-\widetilde{\alpha}_1)\psi''_1+...+(\alpha^i_{p-1}-\widetilde{\alpha}_{p-1})\psi''_{p-1}+\alpha_p^iv^i$.
Or
$\psi'_p=\beta_1^i\psi''_1+...+\beta_{p-1}^i\psi''_{p-1}+\beta^i_pv^i$,
where $\beta_j^i=\alpha^i_j-\widetilde{\alpha}_j$, $j=1,...,p-1$,
$\beta_p^i=\alpha_p^i$,

By lemma 1 we have the following: in unexceptional point $z_i$ the
function $\beta^i_p$ has zero of order $d_i$, and other functions
$\beta^i_j$,  where $j<p$ have zeroes in  $z_i$  of higher order.

\begin{prop} There exist a meromorphic function $h$, such that it has poles only in all,
but  $q'$, unexceptional points $z_{m+1},...,z_t$ with  orders
$d_{m+1},...,d_t$, maybe one additional pole in $P$,  and $q''$
zeroes on $M$ with punctured unexceptional points $z_i$. Here the
number $q'$,$q''$ are such that $q'+q''\leq g$.
\end{prop} \proof

Take a bundle
$F=\xi(z_{m+1})^{d_{m+1}}\otimes...\otimes\xi(z_{t})^{d_{t}}$. It
has a meromorphic section $s$  with a single pole in  $P$  and less
or equal than $g$ zeroes. Also it has a canonical section $s^0$ with
zeroes in unexceptional points $z_{m+1},...,z_t$ of orders
$d_{m+1},...,d_t$. Then $h=\frac{s}{s^0}$ is a meromorphic function.
The point $P$ is its pole.
 If an unexceptional point
$z_i$ isn't a zero of $s$,  than this function has in $z_i$ a pole
of order exactly $d_i$. If an unexceptional point $z_i$ is a zero of
$s$, than $h$ has a pole of lower order (or even a zero). The number
of such points we denote as $q'$.

Let $q''$ be a number of zeroes of $s$ on $M$ with punctured
unexceptional points $z_i$. As $q'+q''$ is a number of zeroes of
$s$, which is less or equal than $g$, we get $q'+q''\leq g$.
\endproof

 Put $\psi''_p=h\psi'_p$. Of course
 $<\psi''_1,...,\psi''_p>=<\psi_1,...,\psi_p>$ and $\psi''_p$ can
 have a pole only in  $P$.

{\bf In the last step} we shall calculate the number of points where
the sections $\psi''_1,...,\psi''_p$ are linearly dependant.

In some neighborhoods of points $z_i$, $i=m+1,...,t$  the equalities
 $\psi''_p=(h\beta_1^i)\psi''_1+...+(h\beta^i_{p-1})\psi''_{p-1}+(h\beta_p^i)v^i$ take place.
 All the functions $(h\beta_j^i)$ are holomorphic in neighborhoods of unexceptional points  $z_{m+1},...,z_t$.  That's
 why the section $\psi''_p$ is holomorphic in these points and can have a pole only in $P$. If an unexceptional point $z_i$ is such that $h$  has a pole of order exactly $d_i$,
 than the function
 $(h\beta_p^i)$ takes nonzero value in $z_i$. That means that
in a neighborhood of such a point $z_i$ we can reconstruct $v^i$
from
 $\psi''_1,...,\psi''_{p-1},\psi''_p$. As in stalks over points from a neighborhood of
 $z_i$ the  sections $\psi''_1,...,\psi''_{p-1},v^i$ form a bases, then the sections $ \psi''_1,...,\psi''_{p-1},\psi''_p$ also form a bases in stalks over points from some neighborhood
 of such $z_i$.
 So in unexceptional points $z_{m+1},...,z_t$, where the function $h$ has a pole of needed orders, the sections  $\psi''_1,...,\psi''_{p-1},\psi''_p$
are independent.

Let's look now, where the sections $
\psi''_1,...,\psi''_{p-1},\psi''_p$ are dependent.  They are
dependent in points $z_1,...,z_m$, where
 the sections  $\psi''_1,...,\psi''_{p-1}$ are dependent. The number of such points is less or equal then $N(p-1,g)$.
 They can be dependent in exceptional points $z_{m+1},...,z_t$, the number of exceptional points is less or equal than  $g$. They can be
  dependent in those $q'$ unexceptional points from $z_{m+1},...,z_t$, where
 the function $h$  doesn't has poles of needed orders.  These are all points from $z_1,...,z_t$, where the sections $ \psi''_1,...,\psi''_{p-1},\psi''_p$ are
 dependent.

 Let's look, what happens on $M\setminus\{z_1,...,z_t\}$. There is a pole of all functions and sections  $P$, and not more than $q''$  zeroes of $h$.
 If a point is not a zero of $h$, than sections $ \psi''_1,...,\psi''_{p-1},\psi''_p$ are
 independent. The number of such points was denoted as $q''$.

 So, finally we get that
$N(p,g)\leq N(p-1,g)+g+q'+q''\leq  N(p-1,g)+2g$.  We know that
$N(1,g)=g$, so the result is $N(p,g)\leq 2pg-g$.

 If we add a pole
$P$, we get that the number of apparent singularities can be made
less or equal than $2pg-g+1$.

Thus, we have a meromorphic trivialization, such that the set $S$ of
poles of trivializing sections and points, where they are finite but
linearly dependent, consists of less or equal than  $2pg-g+1$
points. According to the proposition 1 we construct a system, the
points of $S$ will be the additional singularities. And if
$a_i\notin S$, than $a_i$ will be a fuchsian singulatity.

The main theorem is proofed.

\endproof

\begin{remark}
If we put $p=2$, than we get an estimate $3pg+1$ and in \cite{B} the
estimate is  $3pg-1$, but we in \cite{B} it is suggested, that the
representation is irreducible.
\end{remark}

\begin{remark}
As it can be seen from the first step (see proposition 1 and the
proof of the proposition 2), the solutions of the resulting system
have power-like singularities in the additional singularities. So,
these additional points are regular singular points.
\end{remark}

I'd like to thank R.R Gontsov and I.V. Vyugin for useful
discussions.

\end{document}